\documentclass[10pt]{amsart}
\usepackage[psamsfonts]{amssymb}
\usepackage{amsthm,amscd}
\usepackage{type1cm}
\usepackage{color}

\newtheorem{thm}{Theorem}

\newtheorem{lem}{Lemma}
\newtheorem{cor}{Corollary}
\newtheorem{pro}{Proposition}
\theoremstyle{remark}

\theoremstyle{definition}

\allowdisplaybreaks[4]

\title[Classification of endo-commutative straight algebras]{A classification of 2-dimensional endo-commutative straight  algebras of rank 1 over  a non-trivial field
}

\author[S.-E. Takahasi]{Sin-Ei Takahasi}
\author[K. Shirayanagi]{Kiyoshi Shirayanagi}
\author[M. Tsukada]{Makoto Tsukada}
\address[S.-E. Takahasi]{Laboratory of Mathematics and Games\\ Katsushika 2-371\\ Funabashi\\ Chiba 273-0032\\ Japan}
\email{sin\_ei1@yahoo.co.jp}
\address[K. Shirayanagi, M. Tsukada]{Department of Information Science\\ Toho University\\ Miyama 2-2-1\\ Funabashi\\ Chiba 274-8510\\ Japan}
\email[K. Shirayanagi (Corresponding author)]{kiyoshi.shirayanagi@is.sci.toho-u.ac.jp}
\email[M. Tsukada]{tsukada@is.sci.toho-u.ac.jp}

\makeatletter
\@namedef{subjclassname@2020}{\textup{2020} Mathematics Subject Classification}
\makeatother

\subjclass[2020]{Primary 17A30; Secondary 17D99, 13A99}
\keywords{Nonassociative algebras, Endo-commutative algebras, Commutative algebras, Curled algebras, Straight algebras.}

\begin{document}

\maketitle

\begin{abstract}
An endo-commutative algebra is a nonassociative algebra in which the square mapping preserves multiplication. In this paper, we give a complete classification of 2-dimensional endo-commutative straight algebras of rank one over an arbitrary non-trivial field, where a straight algebra of dimension 2 satisfies the condition that there exists an element $x$ such that $x$ and $x^2$ are linearly independent. We list all multiplication tables of the algebras up to isomorphism.
\end{abstract}

\section{Introduction}\label{sec:intro}
Let $A$ be a nonassociative algebra.  The square mapping $x\mapsto x^2$ from $A$ to itself yields various important concepts of $A$, such as square-rootability, anti-commutativity and zeropotency.  We define $A$ to be {\it{endo-commutative}}, if the square mapping of $A$ preserves multiplication, that is, $x^2y^2=(xy)^2$ holds for all $x, y\in A$.
This terminology comes from the identity $(xx)(yy)=(xy)(xy)$ that depicts the {\it innerly} commutative property.
The concept of endo-commutativity was first introduced in \cite{TST1}, where we gave a complete classification of endo-commutative algebras of dimension 2 over the trivial field $\mathbb F_2$ of two elements.

We separate two-dimensional algebras into two categories: {\it{curled}} and {\it{straight}}.  That is, a 2-dimensional algebra is curled if the square of any element $x$ is a scalar multiple of $x$, otherwise it is straight.  In \cite{TST2} we gave a complete classification of 2-dimensional endo-commutative curled algebras  over an arbitrary non-trivial field.  Despite the difference in canonical representation, we note that a similar result appears in Asrorov-Bekbaev-Rakhimov \cite{ABR}.  Research related to curled algebras can be found in \cite{level2,length1}. 

The aim of this paper is to completely classify 2-dimensional endo-commutative straight algebras of rank 1 over an arbitrary non-trivial field, where the rank of an algebra means the rank of its structure matrix, that is, the matrix of structure constants with respect to a linear base.  The strategy for the classification is based on that of Kobayashi \cite{2dim-comm-char2}.  We can find classifications of associative algebras of dimension 2 over the real and complex number fields in \cite{onkochishin}. For other studies on 2-dimensional algebras, see \cite{moduli,ABR,2dim,variety,classification}. 

Our main theorem states that 2-dimensional endo-commutative straight algebras of rank 1 over an arbitrary non-trivial field $K$ are classified into
three families of specific algebras and one family parameterized by a parameter that takes all values in $K$, up to isomorphism.
The details are described in Theorem 1.  As an application of the main theorem, we further give a complete classification of 2-dimensional endo-commutative straight algebras of rank 1 over an arbitrary non-trivial field in each case of commutative, anti-commutative and associative algebras.  We also show that such algebras are non-unital.  The details are described in Corollary 2.  

\section{Isomorphism criterion for 2-dimensional algebras and characterization of 2-dimensional endo-commutative algebras}\label{sec:iso-criterion}

Let $K$ be a non-trivial field.  For any $X=\small{\begin{pmatrix}a&b\\c&d\end{pmatrix}}\in GL_2(K)$, define
\[
\widetilde{X}=\begin{pmatrix}a^2&b^2&ab&ab\\c^2&d^2&cd&cd\\ac&bd&ad&bc\\ac&bd&bc&ad\end{pmatrix},
\]
where $GL_n(K)$ is the general linear group consisting of nonsingular $n\times n$ matrices over $K$.  Then the following result can be found in \cite{TST2}:

\begin{lem}
The mapping $X\mapsto\widetilde{X}$ is a group-isomorphism from $GL_2(K)$ into $GL_4(K)$ with $|\widetilde X|=|X|^4$.
\end{lem}

Let $A$ be a 2-dimensional algebra over $K$ with a linear base $\{e, f\}$.  Then $A$ is determined by the multiplication table $\begin{pmatrix}e^2&ef\\fe&f^2\end{pmatrix}$.  We write
\[
\left\{
    \begin{array}{@{\,}lll}
     e^2=a_1e+b_1f\\
     f^2=a_2e+b_2f\\
     ef=a_3e+b_3f\\
     fe=a_4e+b_4f\\
   \end{array}
  \right. 
\]
with $a_i, b_i\in K\, \, (1\le i\le4)$ and the matrix
$\begin{pmatrix}a_1&b_1\\a_2&b_2\\a_3&b_3\\a_4&b_4\end{pmatrix}$ is called the structure matrix of $A$ with respect to the base $\{e, f\}$. 

We hereafter will freely use the same symbol $A$ for the matrix and for the algebra because the algebra $A$ is determined by its structure matrix.  Then the following result can be found in \cite{TST2}, and this shows the isomorphism criterion for 2-dimensional algebras:

\begin{pro}
Let $A$ and $A'$ be 2-dimensional algebras over $K$.  Then $A$ and $A'$ are isomorphic iff there is $X\in GL_2(K)$ such that 
\begin{equation}
A'=\widetilde{X^{-1}}AX.
\end{equation}
\end{pro}

\vspace{2mm}

\begin{cor} Let $A$ and $A'$ be 2-dimensional algebras over $K$.  If $A$ and $A'$ are isomorphic, then ${\rm{rank}}\, A={\rm{rank}}\, A'$.
\end{cor}
\vspace{2mm}


For any $a_1, b_1, a_2, b_2, a_3, b_3, a_4, b_4\in K$, we consider the following system of cubic equations:

\begin{equation}
\left\{\begin{array}{@{\,}lll}   
a_1^2a_2+b_1a_2b_2+a_1b_2a_3+b_1a_2a_4=a_1a_3^2+a_2b_3^2+a_3^2b_3+a_3b_3a_4\\
a_1^2a_2+b_1a_2b_2+b_1a_2a_3+a_1b_2a_4=a_1a_4^2+a_2b_4^2+a_3a_4b_4+a_4^2b_4\\
a_1^2a_4+b_1a_4^2+b_1a_2b_4+a_1a_3b_4=a_1^2a_3+b_1a_2b_3+b_1a_3^2+a_1b_3a_4\\
a_2(a_1a_4+a_4b_4+b_2b_4)=a_2(a_1a_3+b_2b_3+a_3b_3)\\
a_1b_1a_2+b_1b_2^2+a_1b_2b_3+b_1a_2b_4=b_1a_3^2+b_2b_3^2+a_3b_3^2+a_3b_3b_4\\
a_1b_1a_2+b_1b_2^2+b_1a_2b_3+a_1b_2b_4= b_1a_4^2+b_2b_4^2+b_3a_4b_4+a_4b_4^2\\
b_1(a_1a_4+a_4b_4+b_2b_4)=b_1(a_1a_3+b_2b_3+a_3b_3)\\
b_1a_2a_4+b_2b_3a_4+b_2^2b_4+a_2b_4^2=b_1a_2a_3+b_2^2b_3+a_2b_3^2+b_2a_3b_4.  
\end{array} \right. 
\end{equation}
The following result can be found in \cite{TST2}, and this gives a necessary and sufficient condition for 2-dimensional algebra $A$ over $K$ to be endo-commutative:
\begin{pro}
Let $A$ be a 2-dimensional algebra over $K$ with structure matrix $\small{\begin{pmatrix}a_1&b_1\\a_2&b_2\\a_3&b_3\\a_4&b_4\end{pmatrix}}$.  Then $A$ is endo-commutative iff  the scalars $a_1, b_1, a_2, b_2, a_3, b_3, a_4, b_4$ satisfy $\rm{(2)}$.
\end{pro}

\section{2-dimensional endo-commutative straight algebras}\label{sec:straight}

Recall that a 2-dimensional algebra is curled if the square of any element $x$ is a scalar multiplication of $x$, otherwise it is straight.  Let $A$ be a 2-dimensional
straight algebra over a non-trivial field $K$ with linear base $\{e, f\}$.  By replacing the bases, we may assume that $e^2=f$.  Write $f^2=pe+qf, ef=ae+bf$ and $fe=ce+df$, where $a, b, c, d, p, q\in K$, and hence the structure matrix of $A$ is 
\begin{equation}
S(p, q, a, b, c, d)\equiv\begin{pmatrix}0&1\\p&q\\a&b\\c&d\end{pmatrix}.
\end{equation}
Conversely, the algebra $S(p, q, a, b, c, d)$ defined by (3) is always straight.  Taking $a_1=0, b_1=1, a_2=p, b_2=q, a_3=a, b_3=b, a_4=c$ and $b_4=d$ in (2),  it becomes
\[
(\sharp_1)\, \left\{\begin{array}{@{\,}lll}   
pq+pc=pb^2+a^2b+abc\\
pq+pa=pd^2+acd+c^2d\\
c^2+pd=pb+a^2\\
p(cd+qd)=p(qb+ab)\\
q^2+pd=a^2+qb^2+ab^2+abd\\
q^2+pb= c^2+qd^2+bcd+cd^2\\
cd+qd=qb+ab\\
pc+qbc+q^2d+pd^2=pa+q^2b+pb^2+qad.  
\end{array} \right. 
\]
But since the seventh equation implies the fourth equation, $(\sharp_1)$ is rewritten as 
\[
(\sharp_2)\, \left\{\begin{array}{@{\,}lll}   
pq+pc=pb^2+a^2b+abc\cdots{\rm{(i_2)}}\\
pq+pa=pd^2+acd+c^2d\cdots{\rm{(ii_2)}}\\
p(d-b)=a^2-c^2\cdots{\rm{(iii_2)}}\\
q^2+pd=a^2+qb^2+ab^2+abd\cdots{\rm{(iv_2)}}\\
q^2+pb= c^2+qd^2+bcd+cd^2\cdots{\rm{(v_2)}}\\
cd+qd=qb+ab\cdots{\rm{(vi_2)}}\\
pc+qbc+q^2d+pd^2=pa+q^2b+pb^2+qad\cdots{\rm{(vii_2)}}.  
\end{array} \right. 
\]
However, with the help of (iii$_2$), we see that (iv$_2$)$-$(v$_2$) yields 
\[
(b+d)\{q(b-d)+ab-cd\}=0.
\]
Then $(\sharp_2)$ is rewritten as
\[
(\sharp_3)\, \left\{\begin{array}{@{\,}lll}   
pq+pc=pb^2+a^2b+abc\cdots{\rm{(i_3)}}\\
pq+pa=pd^2+acd+c^2d\cdots{\rm{(ii_3)}}\\
p(d-b)=a^2-c^2\cdots{\rm{(iii_3)}}\\
q^2+pd=a^2+qb^2+ab^2+abd\cdots{\rm{(iv_3)}}\\
(b+d)\{q(b-d)+ab-cd\}=0\cdots{\rm{(v_3)}}\\
q(d-b)=ab-cd\cdots{\rm{(vi_3)}}\\
pc+qbc+q^2d+pd^2=pa+q^2b+pb^2+qad\cdots{\rm{(vii_3)}}.  
\end{array} \right. 
\]
Note that (vi$_3$) obviously implies (v$_3$).  Then $(\sharp_3)$ is rewritten as
\[
(\sharp_4)\left\{\begin{array}{@{\,}lll}   
pq+pc=pb^2+a^2b+abc\cdots{\rm{(i_4)}}\\
pq+pa=pd^2+acd+c^2d\cdots{\rm{(ii_4)}}\\
p(d-b)=a^2-c^2\cdots{\rm{(iii_4)}}\\
q^2+pd=a^2+qb^2+ab^2+abd\cdots{\rm{(iv_4)}}\\
q(d-b)=ab-cd\cdots{\rm{(vi_4)}}\\
pc+qbc+q^2d+pd^2=pa+q^2b+pb^2+qad\cdots{\rm{(vii_4)}}.  
\end{array} \right. 
\]
By calculating (i$_4$)$-$(ii$_4$), we get from (vi$_4$) that
\begin{align*}
p(c-a)&=p(b^2-d^2)+(a+c)(ab-cd)\\
&=p(b^2-d^2)+(a+c)q(d-b)\, \, ({\rm{by\, \, (vi_4)}})\\
&=(b-d)\{p(b+d)-q(a+c)\}.
\end{align*}
Then $(\sharp_4)$ is rewritten as 
\[
(\sharp_5)\, \left\{\begin{array}{@{\,}lll}   
pq+pc=pb^2+a^2b+abc\cdots{\rm{(i_5)}}\\
p(c-a)=(b-d)\{p(b+d)-q(a+c)\}\cdots{\rm{(ii_5)}}\\
p(d-b)=a^2-c^2\cdots{\rm{(iii_5)}}\\
q^2+pd=a^2+qb^2+ab^2+abd\cdots{\rm{(iv_5)}}\\
q(d-b)=ab-cd\cdots{\rm{(vi_5)}}\\
pc+qbc+q^2d+pd^2=pa+q^2b+pb^2+qad\cdots{\rm{(vii_5)}}.  
\end{array} \right. 
\]
Note that (ii$_5$) and (vi$_5$) imply (vii$_5$).  In fact, we have from (ii$_5$) and (vi$_5$) that
\begin{align*}
p(c-a)&=(b-d)(pb+pd-qa-qc)\, \, (\rm{by\, \, (ii_5)})\\
&=pb^2+pbd-qab-qbc-pbd-pd^2+qad+qcd\\
&=-qab+qcd+pb^2-qbc-pd^2+qad\\
&=-q(ab-cd)+pb^2-qbc-pd^2+qad\\
&=q^2b-q^2d+pb^2-qbc-pd^2+qad\, \, (\rm{by\, \, (vi_5)}),
\end{align*}
which implies (vii$_5$).
Therefore $(\sharp_5)$ is rewritten as 
\begin{equation}
\left\{\begin{array}{@{\,}lll}   
pq+pc=pb^2+a^2b+abc\\
p(c-a)=(b-d)\{p(b+d)-q(a+c)\}\\
p(d-b)=a^2-c^2\\
q^2+pd=a^2+qb^2+ab^2+abd\\
q(d-b)=ab-cd.\\  
\end{array} \right. 
\end{equation}
Namely, if $(a_1, b_1, a_2, b_2, a_3, b_3, a_4, b_4)=(0, 1, p, q, a, b, c, d)$, then (2) is rewritten as (4).  Therefore by Proposition 2, we have the following:

\begin{lem}
Let $p, q, a, b, c, d\in K$.  Then the straight algebra $S(p, q, a, b, c, d)$ is endo-commutative iff the point $(p, q, a, b, c, d)\in K^6$ satisfies {\rm{(4)}}.
\end{lem}

Next, we find a necessary and sufficient condition for straight algebras $S(p, q, a, b, c, d)$ and $S(p', q', a', b', c', d')$ to be isomorphic.  By Proposition 1, we see that $S(p, q, a, b, c, d)\\\cong S(p', q', a', b', c', d')$ holds iff there is $X=\begin{pmatrix}x&y\\z&w\end{pmatrix}\in GL_2(K)$ such that 
\[
\widetilde X\begin{pmatrix}0&1\\p'&q'\\a'&b'\\c'&d'\end{pmatrix}=\begin{pmatrix}0&1\\p&q\\a&b\\c&d\end{pmatrix}X.
\]
Note that
\begin{align*}
\widetilde X\begin{pmatrix}0&1\\p'&q'\\a'&b'\\c'&d'\end{pmatrix}&=\begin{pmatrix}x^2&y^2&xy&xy\\z^2&w^2&zw&zw\\xz&yw&xw&yz\\xz&yw&yz&xw\end{pmatrix}\begin{pmatrix}0&1\\p'&q'\\a'&b'\\c'&d'\end{pmatrix}\\
&=\begin{pmatrix}p'y^2+(a'+c')xy&x^2+q'y^2+(b'+d')xy\\
p'w^2+(a'+c')zw&z^2+q'w^2+(b'+d')zw\\
p'yw+a'xw+c'yz&xz+q'yw+b'xw+d'yz\\
p'yw+a'yz+c'xw&xz+q'yw+b'yz+d'xw\end{pmatrix}\\
\end{align*}
and
\[
\begin{pmatrix}0&1\\p&q\\a&b\\c&d\end{pmatrix}X=\begin{pmatrix}0&1\\p&q\\a&b\\c&d\end{pmatrix}\begin{pmatrix}x&y\\z&w\end{pmatrix}=\begin{pmatrix}z&w\\px+qz&py+qw\\ax+bz&ay+bw\\cx+dz&cy+dw\end{pmatrix}
\]
Then we see that $\widetilde X\begin{pmatrix}0&1\\p'&q'\\a'&b'\\c'&d'\end{pmatrix}=\begin{pmatrix}0&1\\p&q\\a&b\\c&d\end{pmatrix}X$ holds iff the following eight equations hold:
\begin{equation}
\left\{\begin{array}{@{\,}lll}
p'y^2+(a'+c')xy=z\\
x^2+q'y^2+(b'+d')xy=w\\
p'w^2+(a'+c')zw=px+qz\\
z^2+q'w^2+(b'+d')zw=py+qw\\
p'yw+a'xw+c'yz=ax+bz\\
xz+q'yw+b'xw+d'yz=ay+bw\\
p'yw+a'yz+c'xw=cx+dz\\
xz+q'yw+b'yz+d'xw=cy+dw.
\end{array} \right. 
\end{equation}
Therefore we have the following:

\begin{lem}
The algebras $S(p, q, a, b, c, d)$ and $S(p', q', a', b', c', d')$ are isomorphic iff there are $x, y, z, w\in K$ with $xw-yz\ne0$ such that  ${\rm{(5)}}$ holds.
\end{lem}
\vspace{2mm}

\section{Classification of 2-dimensional endo-commutative straight algebras of rank 1}\label{sec:classification}

Denote by $\mathcal{EC}$ the family of all endo-commutative algebras of dimension 2 over a non-trivial field $K$, and define 
\[
\mathcal{ECS}_1=\{S(p, q, a, b, c, d)\in\mathcal{EC} : {\rm{rank}}\, S(p, q, a, b, c, d)=1\}.
\]
Since the algebra $S(p, q, a, b, c, d)$ has rank one iff $p=a=c=0$, it follows from Lemma 2 that
\[
\mathcal{ECS}_1=\{S(0, q, 0, b, 0, d) : (0, q, 0, b, 0, d)\in K^6\, \, {\rm{satisfies\, \, (4)}}\}.
\]
Put
\[
\mathcal{ECS}_{10}=\{S(0, 0, 0, b, 0, d) : (0, 0, 0, b, 0, d)\in K^6\, \, {\rm{satisfies\, \, (4)}}\}
\]
and
\[
\mathcal{ECS}_{11}=\{S(0, q, 0, b, 0, d) : (0, q, 0, b, 0, d)\in K^6\, \, {\rm{satisfies\, \, (4)}}\, \, {\rm{and}}\, \, q\ne0\}.
\]
Then we have the disjoint union:
\[
\mathcal{ECS}_1=\mathcal{ECS}_{10}\sqcup\mathcal{ECS}_{11}.
\]
Now we easily see that $(0, q, 0, b, 0, d)\in K^6$ satisfies (4) iff $\left\{\begin{array}{@{\,}lll}   
q^2=qb^2\\
q(d-b)=0. 
\end{array} \right.$ Then we have
\[
\mathcal{ECS}_{10}=\{S(0, 0, 0, b, 0, d) : b, d\in K\}
\]
and
\[
\mathcal{ECS}_{11}=\{S(0, b^2, 0, b, 0, b) : b\in K\, \, {\rm{and}}\, \, b\ne0 \}.
\]
By Lemma 3, we see that $S(0, q, 0, b, 0, d)\cong S(0, q', 0, b', 0, d')$ iff there are $x, y, z, w\in K$ with $xw\ne  yz$ such that
\begin{equation}
\left\{\begin{array}{@{\,}lll}
z=0\\
x^2+q'y^2+(b'+d')xy=w\\
q'w^2=qw\\
q'yw+b'xw=bw\\
q'yw+d'xw=dw.
\end{array} \right. 
\end{equation}
If $z=0$ and $xw\ne yz$, then $xw\ne0$.  Then (6) is rewritten as
\[
\left\{\begin{array}{@{\,}lll}
z=0\\
x^2+q'y^2+(b'+d')xy=w\\
q'w=q\\
q'y+b'x=b\\
q'y+d'x=d.
\end{array} \right. 
\]
Consequently, we have the following:

\begin{lem}
The algebras $S(0, q, 0, b, 0, d)$ and $S(0, q', 0, b', 0, d')$ are isomorphic iff there are $x, y, w\in K$ with $xw\ne0$ such that 
\begin{equation}
\left\{\begin{array}{@{\,}lll}
x^2+q'y^2+(b'+d')xy=w\\
q'w=q\\
q'y+b'x=b\\
q'y+d'x=d
\end{array} \right. 
\end{equation}
holds.
\end{lem}

Looking at the second equation in (7),  we see from Lemma 4 that each algebra in $\mathcal{ECS}_{10}$ is not isomorphic to any algebra in $\mathcal{ECS}_{11}$.
Therefore, we will first separately classify the algebras in $\mathcal{ECS}_{10}$ and the algebras in $\mathcal{ECS}_{11}$, and then
combine them into the desired classification of $\mathcal{ECS}_1$.

\vspace{3mm}

\noindent (I) Classification of $\mathcal{ECS}_{10}$.

\vspace{3mm}

First we have the following:

\begin{lem}
  Any algebra in $\mathcal{ECS}_{10}$ is isomorphic to one of the algebras 
  $S(0, 0, 0, 0, 0, 0)$, $S(0, 0, 0, 0, 0, 1)$ and $S(0, 0, 0, 1, 0, \lambda)\, \, (\lambda\in K)$.
\end{lem}

\begin{proof}

Let $S(0, 0, 0, b, 0, d), S(0, 0, 0, b', 0, d')\in\mathcal{ECS}_{10}$ be arbitrary.  Then we see from Lemma 4 that $S(0, 0, 0, b, 0, d)\cong S(0, 0, 0, b', 0, d')$ iff there are $x, y, w\in K$ with $xw\ne0$ such that 
\begin{equation}
\left\{\begin{array}{@{\,}lll}
x^2+(b'+d')xy=w\\
b'x=b\\
d'x=d.
\end{array} \right. 
\end{equation}

(i) The case where $b=d=0$.  Trivially, $S(0, 0, 0, b, 0, d)=S(0, 0, 0, 0, 0, 0)$.
\vspace{2mm}

(ii) The case where $b=0$ and $d\ne0$.  Put $x=d, y=0$ and $w=d^2$.  Then $xw=d^3\ne0$ and (8) holds with $b'=0$ and $d'=1$.  Then we have $S(0, 0, 0, b, 0, d)\cong S(0, 0, 0, 0, 0, 1)$.
\vspace{2mm}

(iii) The case where $b\ne0$.   Put $x=b, y=0$ and $w=b^2$.  Then $xw=b^3\ne0$ and (8) holds with $b'=1$ and $d'=d/b$.  Then we have $S(0, 0, 0, b, 0, d)\cong S(0, 0, 0, 1, 0, d/b)$.

Therefore the desired result follows immediately from (i), (ii) and (iii) above.
\end{proof}

We next investigate the relationship of isomorphism between these algebras.  Then we have the following:

\begin{pro}
  The algebras in $\mathcal{ECS}_{10}$ are classified into three families $S(0, 0, 0, 0, 0, 0)$, $S(0, 0, 0, 0, 0, 1)$ and
  $S(0, 0, 0, 1, 0, \lambda)\, (\lambda\in K)$ up to isomorphism.
\end{pro}

\begin{proof}
(i) $S(0, 0, 0, 0, 0, 0)\ncong S(0, 0, 0, 0, 0, 1)$.  Suppose on the contrary that \\$S(0, 0, 0, 0, 0, 0)\cong S(0, 0, 0, 0, 0, 1)$.  Then by Lemma 4, there are $x, y, w\in K$ with $xw\ne0$ such that $\left\{\begin{array}{@{\,}lll}
x^2+xy=w\\
x=0.
\end{array} \right. $  These equations immediately imply $x=w=0$, which contradicts $xw\ne0$.
\vspace{2mm}

(ii) $S(0, 0, 0, 0, 0, 0)\ncong S(0, 0, 0, 1, 0, \lambda)$, $\lambda\in K$.  Suppose on the contrary that $S(0, 0, 0, 0, 0, 0)\cong S(0, 0, 0, 1, 0, \lambda)$.  Then by Lemma 4, there are $x, y, w\in K$ with $xw\ne0$ such that $\left\{\begin{array}{@{\,}lll}
x^2+(1+\lambda)xy=w\\
x=0\\
\lambda x=0.
\end{array} \right.$  The first two equations imply $x=w=0$, which contradicts $xw\ne0$.
\vspace{2mm}

(iii) $S(0, 0, 0, 0, 0, 1)\ncong S(0, 0, 0, 1, 0, \lambda)$, $\lambda\in K$.  Suppose on the contrary that $S(0, 0, 0, 0, 0, 1)\cong S(0, 0, 0, 1, 0, \lambda)$.  Then by Lemma 4, there are $x, y, w\in K$ with $xw\ne0$ such that $\left\{\begin{array}{@{\,}lll}
x^2+(1+\lambda)xy=w\\
x=0\\
\lambda x=1.
\end{array} \right.$  However, the second and third equations contradict each other.
\vspace{2mm}

(iv) $S(0, 0, 0, 1, 0, \lambda)\ncong S(0, 0, 0, 1, 0, \lambda')$ with $\lambda\neq\lambda'$. Suppose on the contrary that $S(0, 0, 0, 1, 0, \lambda)\cong S(0, 0, 0, 1, 0, \lambda')$.  Then by Lemma 4, there are $x, y, w\in K$ with $xw\ne0$ such that $\left\{\begin{array}{@{\,}lll}
x^2+(1+\lambda')xy=w\\
x=1\\
\lambda'x=\lambda.
\end{array} \right. $  But the second and third equations imply $\lambda=\lambda'$, a contradiction.

By (i), (ii), (iii), (iv) and Lemma 5, we obtain the desired result.
\end{proof}

\noindent (II) Classification of $\mathcal{ECS}_{11}$.

\vspace{3mm}

We have the following:

\begin{pro}
Any algebra in $\mathcal{ECS}_{11}$ is isomorphic to $S(0, 1, 0, 1, 0, 1)$.
\end{pro}

\begin{proof}
Let $b, b'\ne0$ and $S(0, b^2, 0, b, 0, b), S(0, b'^2, 0, b', 0, b')\in\mathcal{ECS}_{11}$.  Then we see from Lemma 4 that $S(0, b^2, 0, b, 0, b)\cong S(0, b'^2, 0, b', 0, b')$ iff there are $x, y, w\in K$ with $xw\ne0$ such that
\begin{equation}
\left\{\begin{array}{@{\,}lll}
x^2+b'^2y^2+2b'xy=w\\
b'^2w=b^2\\
b'^2y+b'x=b.
\end{array} \right. 
\end{equation}
Put $x=b, y=0$ and $w=b^2$.  Then $xw=b^3\ne0$ and (9) holds with $b'=1$, hence $S(0, b^2, 0, b, 0, b)\cong S(0, 1, 0, 1, 0, 1)$ from the above argument.
\end{proof} 
\vspace{3mm}

Recall that  each algebra in $\mathcal{ECS}_{10}$ is not isomorphic to any algebra in $\mathcal{ECS}_{11}$.  Then, by Propositions 3 and 4, we have the following:

\begin{thm}
Up to isomorphism, 2-dimensional endo-commutative straight algebras of rank 1 over a non-trivial field $K$ are classified into 4 families
\[
 S(0, 0, 0, 0, 0, 0), S(0, 0, 0, 0, 0, 1), \{S(0, 0, 0, 1, 0, \lambda)\}_{ \lambda\in K}\, \, {\rm{and}}\, \, S(0, 1, 0, 1, 0, 1)
\]
with multiplication tables on a linear base $\{e, f\}$ defined by
\[
\begin{pmatrix}f&0\\0&0\end{pmatrix}, \begin{pmatrix}f&0\\f&0\end{pmatrix}, \begin{pmatrix}f&f\\\lambda f&0\end{pmatrix}\, \, {\rm{and}}\, \, \begin{pmatrix}f&f\\f&f\end{pmatrix},
\]
respectively. 
\end{thm}

\begin{cor}
Suppose that $A$ is a 2-dimensional endo-commutative straight algebra of rank 1 over a non-trivial field $K$. Then 
\vspace{2mm}

{\rm{(i)}} $A$ is non-unital. 
\vspace{2mm}

{\rm{(ii)}} $A$ is commutative iff $A\cong S(0, 0, 0, 0, 0, 0), A\cong S(0, 0, 0, 1, 0, 1)$ or $A\cong S(0, 1, 0, 1, 0, 1)$. 
\vspace{2mm}

{\rm{(iii)}} If ${\rm{char}} K\ne2$, then $A$ is not anti-commutative. 
\vspace{2mm}

{\rm{(iv)}} $A$ is associative iff $A\cong S(0, 0, 0, 0, 0, 0)$ or $A\cong S(0, 1, 0, 1, 0, 1)$. 
\end{cor}

\begin{proof}
Let $\{e, f\}$ be a linear base of $A$ with $e^2=f$.
\vspace{2mm}

(i) Assume that $A$ has an identity element $u=\alpha e+\beta f$, where $\alpha, \beta\in K$.  If $A$ has the multiplication table $\begin{pmatrix}f&0\\0&0\end{pmatrix}$, then $e=ue=\alpha e^2+\beta fe=\alpha f$, a contradiction.  If  $A$ has the multiplication table $\begin{pmatrix}f&0\\f&0\end{pmatrix}$ or  $\begin{pmatrix}f&f\\f&f\end{pmatrix}$, then  $e=ue=\alpha e^2+\beta fe=(\alpha+\beta)f$, a contradiction. Also if $A$ has the multiplication table $\begin{pmatrix}f&f\\\lambda f&0\end{pmatrix}$, then $e=ue=\alpha e^2+\beta fe=(\alpha+\beta\lambda)f$, a contradiction. Therefore assertion (i) follows from Theorem 1.
\vspace{2mm}

(ii) Since $A$ is commutative iff the multiplication table of $A$ is symmetric, assertion (ii) follows from Theorem 1.
\vspace{2mm}

(iii) A simple consideration shows that $A$ is anti-commutative iff 
\begin{equation}
2acf+(bc+ad)(ef+fe)+2bd\,f^2=0
\end{equation}
holds for all $a, b, c, d\in K$. 

Suppose ${\rm{char}} K\ne2$ and $A$ is anti-commutative.   If $A$ has the multiplication table $\begin{pmatrix}f&0\\0&0\end{pmatrix}$, then (10) implies that $2acf=0$ for all $a, c\in K$, a contradiction.  If $A$ has the multiplication table $\begin{pmatrix}f&0\\f&0\end{pmatrix}$, then (10) implies that $(2ac+bc+ad)f=0$  for all $a, b, c, d\in K$, a contradiction.  If $A$ has the multiplication table $\begin{pmatrix}f&f\\\lambda f&0\end{pmatrix}$, then (10) implies that $\{2ac+(1+\lambda)(bc+ad)\}f=0$ for all $a, b, c, d\in K$, a contradiction.  Also if f $A$ has the multiplication table $\begin{pmatrix}f&f\\f&f\end{pmatrix}$, then (10) implies that $2(a+b)(c+d)f=0$ for all $a, b, c, d\in K$, a contradiction.  Therefore assertion (iii) follows from Theorem 1.
\vspace{2mm}

(iv)  A simple consideration shows that $A$ is associative iff 
\begin{equation}
\left\{
    \begin{array}{@{\,}lll}
     e^2e=ee^2\\
     (ef)e=e(fe)\\
     (fe)e=fe^2\\
     f^2e=f(fe)\\
     e^2f=e(ef)\\
     (ef)f=ef^2\\
     (fe)f=f(ef)\\
     f^2f=ff^2
   \end{array}
  \right. 
\end{equation}
holds.   If $A$ has the multiplication table $\begin{pmatrix}f&0\\0&0\end{pmatrix}$, then (11) holds since both sides of each equation in (11) are 0,
  and so $A$ is associative.  If $A$ has the multiplication table $\begin{pmatrix}f&0\\f&0\end{pmatrix}$, then the third equation of (11) implies that $f=fe=(fe)e=fe^2=f^2=0$, a contradiction, hence $A$ is not associative.  If $A$ has the multiplication table $\begin{pmatrix}f&f\\\lambda f&0\end{pmatrix}$, then the fifth equation of (11) implies that $f=ef=e(ef)=e^2f=f^2=0$, a contradiction, hence $A$ is not associative.  Finally, if $A$ has the multiplication table $\begin{pmatrix}f&f\\f&f\end{pmatrix}$,
then (11) holds since both sides of each equation in (11) are $f$, hence $A$ is associative.  Summarizing the above, we obtain assertion (iv).
\end{proof}

\noindent Remark. Among the 2-dimensional endo-commutative straight algebras of rank 1, Corollary 2 implies that only $S(0, 0, 0, 0, 0, 1)$ and $S(0, 0, 0, 1, 0, \lambda)\, \, (\lambda\ne1)$ are {\it purely} endo-commutative in the sense that they are neither commutative nor anti-commutative nor associative.


\vspace{2mm}
\begin{thebibliography}{99}

\bibitem{moduli}
A. Ananin and A. Mironov, The moduli space of 2-dimensional algebras, Comm. Algebra, \textbf{28}(9) (2000), 4481-4488. 

\bibitem{ABR} D. Asrorov, U. Bekbaev and I. Rakhimov, A complete classification of two-dimensional endo-commutative algebras over an arbitrary field,
  arXiv:2304.00491v1 [math.RA] 2 Apr 2023.

\bibitem{2dim}
M. Goze and E. Remm, 2-dimensional algebras, Afr. J. Math. Phys., \textbf{10}(1) (2011), 81-91. 

\bibitem{variety}
I. Kaygorodov and Y. Volkov, The variety of 2-dimensional algebras over an algebraically closed field, Canad. J. Math., \textbf{71}(4) (2019), 819-842. 

\bibitem{level2}
I. Kaygorodov and Y. Volkov, Complete classification of algebras of level two, Mosc. Math. J., \textbf{19}(3) (2019), 485-521. 

\bibitem{2dim-comm-char2} 
Y. Kobayashi, Characterization of two-dimensional commutative algebras over a field of characteristic 2, preprint (2016).

\bibitem{onkochishin} 
  Y. Kobayashi, K. Shirayanagi, M. Tsukada and S.-E. Takahasi, A complete classification of three-dimensional algebras over $\mathbb R$ and $\mathbb C$
  -- {\it OnkoChishin} (visiting old, learn new), Asian-Eur. J. Math., \textbf{14}(8) (2021), Article ID 2150131 (25 pages). 

\bibitem{length1}
O. Markova, C. Mart\'{i}nez and R. Rodrigues, Algebras of length one, J. Pure Appl. Algebra, \textbf{226}(7) (2022), Article ID 106993, 16 p. 

\bibitem{classification}
H. Petersson, The classification of two-dimensional nonassociative algebras, Results Math., \textbf{37}(1-2) (2000), 120-154. 

\bibitem{TST1} S.-E. Takahasi, K. Shirayanagi, M. Tsukada, A classification of two-dimensional endo-commutative algebras over $\mathbb F_2$, arXiv:2211.04015v1 [math.RA] 8 Nov 2022. 


\bibitem{TST2} S.-E. Takahasi, K. Shirayanagi, M. Tsukada, A classification of endo-commutative curled algebras of dimension 2 over a non-trivial field, arXiv:2304.2510v1 [math.RA] 25 Apr 2023. 

\end{thebibliography}
\end{document}